\documentclass{CSML}
\pdfoutput=1

\usepackage{lastpage}
\newcommand{\dOi}{13(2:11)2017}
\lmcsheading{}{1--\pageref{LastPage}}{}{}%
{Nov.~02, 2016}{Jun.~23, 2017}{}

\usepackage{latexsym, url}
\usepackage{amsthm}
\usepackage{amsmath}
\usepackage{amsfonts}
\usepackage{amssymb}
\usepackage{hyperref}
\hypersetup{hidelinks}
\usepackage[all]{xypic}

\newcommand\ksst{\preceq_\ck}

\newcommand\lsk{LS(\ck)}
\newcommand\lsdk{LS^d(\ck)}
\newcommand\sgr{\triangleright}
\newcommand\sgeq{\trianglerighteq}

\newcommand\op{\operatorname{op}}

\newcommand\Set{\operatorname{\bf Set}}
\newcommand\Met{\operatorname{\bf Met}}

\newcommand\FrAb{\operatorname{\bf FrAb}}
\newcommand\Ab{\operatorname{\bf Ab}}

\newcommand\Presll{\operatorname{{\bf Pres}_\lambda(\cl)}}
\newcommand\Preskk{\operatorname{{\bf Pres}_\kappa(\ck)}}
\newcommand\Preslk{\operatorname{{\bf Pres}_\lambda(\ck)}}

\newcommand\gatp{\operatorname{ga-tp}}

\newcommand\colim{\operatorname{colim}}

\newcommand\monst{\mathfrak {C}}
\newcommand\reals{\mathbb {R}}
\newcommand\ca{\mathcal {A}}

\newcommand\cd{\mathcal {D}}

\newcommand\ck{\mathcal {K}}
\newcommand\ckpa{\mathcal {K^{\bullet\,\bullet}}}
\newcommand\cl{\mathcal {L}}

\begin{document}
\title[Hanf numbers via accessible images]
{Hanf numbers via accessible images}
\author[M. Lieberman and J. Rosick\'{y}]
{M. Lieberman \and J. Rosick\'y}
\address{Department of Mathematics and Statistics,
Masaryk University, 
Czech Republic}
\email{\{lieberman,rosicky\}@math.muni.cz}
\thanks{Supported by the Grant Agency of the Czech Republic under the grant 
               P201/12/G028.} 
\dedicatory{Dedicated to Ji\v r\'i Ad\'amek on the occasion of his 70th birthday.}
 
 
\begin{abstract}
We present several new model-theoretic applications of the fact that, under the assumption that there exists a proper class of almost strongly compact cardinals, the powerful image of any accessible functor is accessible.  In particular, we generalize to the context of accessible categories the recent Hanf number computations of Baldwin and Boney, namely that in an abstract elementary class (AEC) if the joint embedding and amalgamation properties hold for models of size up to a sufficiently large cardinal, then they hold for models of arbitrary size.  Moreover, we prove that, under the above-mentioned large cardinal assumption, every metric AEC is strongly $d$-tame, strengthening a result of Boney and Zambrano and pointing the way to further generalizations.
\end{abstract} 

\maketitle

\section{Introduction}

The connections of accessible categories to problems in abstract model theory, and particularly to those involving abstract elementary classes (AECs) and metric abstract elementary classes (mAECs), have by now been the subject of considerable investigation, e.g. \cite{L}, \cite{BR}, \cite{LR}, \cite{muaecs}, and \cite{LRmetr}.  Accessible categories serve as a common generalization, certainly---AECs, mAECs, and further generalizations including $\mu$-CAECs and $\mu$-AECs,  all amount to accessible categories with additional structure---and offer a perspective in which certain powerful tools for analysing the structure of such classes become evident.  As in \cite{LR}, all of the machinery at work in this paper can be traced back, ultimately, to Theorem 5.5.1 in \cite{MP}: assuming arbitrarily large strongly compact cardinals, given any accessible functor $F:\ck\to\cl$, the powerful image of $F$ (that is, the closure of its image under subobjects in $\cl$) is an accessible category.

The dependence on set theory, while perhaps surprising, is easily motivated.  While the category of abelian groups, $\Ab$, is finitely accessible, the full subcategory of free abelian groups, $\FrAb$, which is the powerful image of the accessible free abelian group functor $F:\Set\to\Ab$, need not be.  In particular, it is not accessible assuming $V=L$; if there are arbitrarily large strongly compact cardinals, on the other hand, it is accessible (\cite{EM}).  The aforementioned theorem of \cite{MP}, which generalizes this fact, admits a further generalization in \cite{BTR}: if $F:\ck\to\cl$ is a $\lambda$-accessible functor that preserves $\mu_\cl$-presentable objects---we recall the definition of $\mu_\cl$ in Section~\ref{sprops} below---and $\kappa$ is a $\mu_\cl$-strongly compact cardinal, the powerful image of $F$ is $\mu_\cl$-accessible.

The applications of the accessibility of powerful images in abstract model theory, both past (the proof in \cite{LR}, via \cite{MP}, of Boney's Theorem on the tameness of AECs under the assumption of arbitrarily large strongly compact cardinals in \cite{B}---subsequently completed in \cite{BU} 4.14 to a grand equivalence of tameness, accessibility of powerful images, and arbitrarily large almost strongly compact cardinals) and present (the arguments to follow, which lean on the theorem of \cite{BTR} above) are united by a single, very simple idea.  In particular, many of the essential properties of abstract classes of structures $\ck$---including the amalgamation and joint embedding properties (abbreviated as AP and JEP, respectively), and the equality of, or distance between, Galois types---reduce to questions about whether diagrams of a particular shape $\ca$ can be extended to diagrams of shape $\ca'$.  In the case of the joint embedding and amalgamation properties, and to an extent their disjoint versions (DJEP, DAP), $\ca$ and $\ca'$ are simple finite categories and the diagrams live entirely in $\ck$---which we now take to be a general accessible category, meaning that we can identify the respective categories of diagrams with the functor categories $\ck^\ca$ and $\ck^{\ca'}$.  This sets up the first, and perhaps clearest, application of the accessibility of powerful images.

The functor whose image we will consider is simply the forgetful functor
\[F:\ck^{\ca'}\to\ck^\ca\]
	which forgets the details of the extension.  Clearly, the category of extendable diagrams of shape $\ca$ corresponds to its (powerful) image.  When this functor is accessible, as it will be in the cases considered below, we know that the image is $\kappa$-accessible for any sufficiently compact cardinal $\kappa$.  If we know that any $\ca$-diagram of $\kappa$-presentable objects is extendible (corresponding to, e.g. the $<\kappa$-JEP), the latter fact---in particular, the closure of the image under $\kappa$-directed colimits---is precisely what we need to infer that every $\ca$-diagram is extendible (corresponding to, e.g. the full JEP).  In Sections~\ref{sprops} and \ref{disj}, we will use this template to generalize Theorem 1.1 of \cite{BaBo}, proving that if one of the JEP, AP, DJEP, or DAP holds in an accessible category $\ck$ for objects of presentability rank up to a sufficiently large cardinal $\kappa$---a Hanf number, of sorts---it must hold for objects of arbitrary presentability rank.

In Section~\ref{mtame}, we turn to the question of $d$-tameness of Galois types in mAECs.  Here again, the problem reduces to the accessibility of a particular category of completable diagrams, and is a straightforward (albeit decidedly more technical) generalization of the corresponding category-theoretic proof for the discrete case, namely conventional AECs (\cite{LR} 5.2).  The discrete case is instructive: in an AEC $\ck$, a Galois type over $M\in\ck$ is an equivalence class of pairs $(f:M\to N, a\in UN)$, where $UN$ denotes the underlying set of $N$, and where pairs $(f_i:M\to N_i, a_i\in UN_i)$ for $i=1,2$ are equivalent if there is an amalgam in which $a_0$ and $a_1$ are identified.  This is, essentially, a pointed version of the picture for the AP, with the category of equivalent pairs as the image of the obvious forgetful functor.  This (powerful) image will be $\kappa$-accessible for sufficiently compact $\kappa$, meaning that the category of equivalent pairs is determined by the subcategory of equivalent pairs over models of size no greater than $\kappa$. In the terminology of AECs, this means that $\ck$ is $\kappa$-tame---see \cite{LR} for details.

In an mAEC $\ck$, the objects have not just underlying sets, but underlying complete metric spaces.  Provided $\ck$ is reasonably well-behaved, we can speak meaningfully of the distance between Galois types $(f_i:M\to N_i, a_i\in UN_i)$ for $i=1,2$, here defined not in an ambient monster model, as is customary, but as the infimum of the distances between the images of $a_0$ and $a_1$ in all possible amalgams of $N_0\leftarrow M\to N_1$ in $\ck$.  For each $\epsilon>0$, we then have a forgetful functor $G_\epsilon$ whose (powerful) image is the category of pairs within distance $\epsilon$ of one another.  By choosing sufficiently compact $\kappa$, we may assume that this image is $\kappa$-accessible, and can thereby force a very strong version of metric tameness (i.e. $d$-tameness).  Namely, we show that for any $\epsilon>0$, a pair of types over an arbitrary model $M\in\ck$ are at distance at least $\epsilon$ from one another only if there is a submodel $N$ of $M$ of density character at most $\kappa$ such that their restrictions to $N$ are already at distance at least $\epsilon$.  This strengthens, at least superficially, the tameness portion of Theorem 3.11 in \cite{BoZ}, which is concerned with the standard $\epsilon-\delta$ definition of $d$-tameness; we prove that one may take $\delta=\epsilon$.  More significantly, this approach points a way to further generalizations of tameness phenomena.

We conclude with a discussion of possible extensions of these results.  In ongoing work with Zambrano, the authors are investigating natural generalizations from metric structures to quantale-valued structures, with the enticing prospect of meaningful applications both in fuzzy, probabilistic contexts, and in the more rarefied realm of model theory over sheaves.

We thank the anonymous referees for their insightful comments on this paper, which have led to noteworthy improvements in the presentation.

\section{Preliminaries}\label{prelim}

We here restate a few basic definitions and results concerning accessible categories and their connections with AECs and mAECs; exhaustive treatments of the former topic can be found in \cite{MP} and \cite{AR}, while the current state of knowledge on the latter is summarized in \cite{LRmetr}.  Additional background on AECs, and the broader context of abstract model theory, can be found in, e.g. \cite{Ba}.

For $\lambda$ a regular cardinal, a {\em $\lambda$-accessible category} is a category $\ck$ with $\lambda$-directed colimits in which every object can be expressed as a $\lambda$-directed colimit of $\lambda$-presentable objects, of which $\ck$ contains only set many, up to isomorphism.  Here an object $N$ in $\ck$ is {\em $\lambda$-presentable} if the hom-functor $\ck(N,-):\ck\to\Set$ preserves $\lambda$-directed colimits; the {\em presentability rank} of an object $N$ in $\ck$, denoted $\pi(N)$, is the least $\lambda$ such that $N$ is $\lambda$-presentable.  Presentability rank serves as a notion of size in an abstract category $\ck$, and connects nicely to the natural notions of size in familiar model-theoretic contexts.

\begin{exas}
(1) Let $\ck$ be an AEC with L\"owenheim-Skolem number $\lsk$, and let $U:\ck\to\Set$ be the underlying set functor.  Then $\ck$ is $\lsk^+$-accessible with directed colimits.  Moreover, for any regular $\kappa>\lsk$, a model $N\in\ck$ is $\kappa$-presentable if and only if $|UN|<\kappa$ (\cite{L} 4.1, 4.3).  In particular, $\pi(N)=|UN|^+$.

(2) Let $\ck$ be an mAEC with L\"owenheim-Skolem number $\lsdk$, and let $U:\ck\to\Met$ be the functor assigning to each model its underlying complete metric space: recall that $\Met$ is the category of complete metric spaces and isometric embeddings. Here the relevant notion of size for models $M\in\ck$ is not the cardinality of $UM$, but its density character, which we denote by $\hbox{dc}(UM)$---this is the reason for the superscript in $\lsdk$.  By \cite{LRmetr} 3.1, $\ck$ is $\lsdk^+$-accessible with directed colimits.  Moreover, for regular $\kappa>\lsdk$, a model $N\in\ck$ is $\kappa$-presentable if and only if $\hbox{dc}(UN)<\kappa$.  In particular, $\pi(N)=\hbox{dc}(UN)^+$.  
\end{exas}

In the examples above, and, indeed, in any accessible category with directed colimits, the categories $\ck$ are not merely accessible in $\lsk$ or $\lsdk$, but in all larger regular cardinals, i.e. they are {\em well-accessible} (\cite{BR} 4.1).  This is not true in general: for $\mu>\lambda$, a $\lambda$-accessible category is $\mu$-accessible if and only if $\mu$ is {\em sharply greater than} $\lambda$:

\begin{defi} For regular cardinals $\lambda$ and $\mu$, we say that $\mu$ is {\em sharply greater than} $\lambda$, denoted $\mu\sgr\lambda$, if $\mu>\lambda$ and for any $\kappa<\mu$ the set $[\kappa]^{<\lambda}$ consisting of subsets of $\kappa$ of size less than $\lambda$, ordered by subset inclusion, contains a cofinal set of size less than $\mu$.	
\end{defi}

A number of equivalent characterizations of this relation can be given, e.g. in \cite{MP} 2.3.  We note that, barring the assumption of GCH, this condition on $\mu$ and $\lambda$ is weaker than the cardinal inequality $\mu^\lambda=\mu$.  The latter condition is not a bad intuition to keep around, however: following \cite{LRmetr} 4.11, if $\mu>2^\lambda$, then $\mu^+\sgr\lambda^+$ if and only if $\mu^\lambda=\mu$.

\begin{defi} Let $F:\ck\to\cl$ be a functor.
\begin{enumerate}
	\item We say that $F$ is {\em $\lambda$-accessible} if it preserves $\lambda$-directed colimits and the categories $\ck$ and $\cl$ are both $\lambda$-accessible.
	\item The {\em powerful image} of $F$ is the smallest subcategory of $\cl$ that contains the image of $F$ and is closed under subobjects.
\end{enumerate} 
\end{defi}

Finally, we recall the definitions of the large cardinals with which we will be working:

\begin{defi} Given cardinals $\kappa\geq\mu$, we say that $\kappa$ is {\em $\mu$-strongly compact} if any $\kappa$-complete filter can be extended to a $\mu$-complete ultrafilter.  A cardinal $\kappa$ is {\em almost strongly compact} if it is $\mu$-strongly compact for all $\mu<\kappa$.\end{defi}

We note that the assumption of $\mu$-strongly compact cardinals for arbitrary $\mu$, preferred in \cite{BTR}, is equivalent to the assumption, preferred in \cite{BU}, of a proper class of almost strongly compact cardinals (see \cite{BTR}, Proposition~2.4).  We note, too, that whereas \cite{BTR} derives the accessibility of powerful images from these equivalent large cardinal assumptions, \cite{BU} shows that this assumption is optimal: the powerful image of any accessible functor is accessible if and only if there is a proper class of almost strongly compact cardinals.

The key result in the argument that follows is Theorem 3.4 in \cite{BTR}, which guarantees that for any well-behaved $\lambda$-accessible functor $F:\ck\to\cl$, its powerful image is $\kappa$-accessible for $\mu_\cl$-strongly compact $\kappa$, where $\mu_\cl$ is a cardinal that depends on both $\lambda$ and $\cl$.  We are now in a position to describe this $\mu_\cl$, following the notation of \cite{BTR} 3.1.  

Let $\Presll$ denote a skeleton of the full subcategory of $\lambda$-presentable objects in $\cl$; that is, a subcategory containing a single object for each isomorphism class of $\lambda$-presentable objects.  Let $\beta=|\Presll|$.  Let $\gamma_\cl$ be the smallest cardinal such that $\gamma_\cl\geq\beta$ and $\gamma_\cl\sgeq\lambda$.  We note that this $\gamma_\cl$ will not be too large: if $\lambda<\beta$, 
\[\gamma_\cl\leq(2^\beta)^+,\]
and clearly $\gamma_\cl=\lambda$ otherwise (\cite{AR} 2.13(3)).  In case $\ck$ is well $\lambda$-accessible, we may simply take $\gamma_\cl=\max(\beta,\lambda)$, noting the parallel between this cardinal and the corresponding bound for AECs used in \cite{BaBo}, i.e. $I(\ck,\lambda)+\lambda$.  Finally, we define $\mu_\cl=(\gamma_\cl^{<\gamma_\cl})^+$.  By design, $\mu_\cl\sgeq\lambda$ (\cite{AR} 2.13(5)). This notation is slightly ambiguous, given that $\cl$ will be $\lambda$-accessible for many $\lambda$, but the choice of $\lambda$ will usually be clear.  In the one place where there is a risk of confusion--- Theorems~\ref{disjoint3} and \ref{disjoint4}---we specify the cardinal parameter in the subscript, e.g. $\mu_{\cl,\lambda}$. 

For the sake of completeness, we state the theorem of \cite{BTR} on which we will depend in the sequel:

\begin{thm}[\cite{BTR} 3.4]\label{btrthm} Let $\lambda$ be a regular cardinal and $\cl$ an accessible category such that there exists a $\mu_\cl$-strongly compact cardinal $\kappa$.  Then the powerful image of any $\lambda$-accessible functor to $\cl$ that preserves $\mu_\cl$-presentable objects is accessible.\end{thm}

We need a slight refinement, however.  In the proof of this theorem in \cite{BTR}, the only obstacle to the powerful image being $\kappa$-accessible is that it may not be the case that $\kappa\sgeq\lambda$.  If, as in the theorem above, one were not concerned with the index of accessibility of the powerful image, one could simply pass to $\kappa'>\kappa$ (hence also $\mu_\cl$-strongly compact) such that $\kappa'\sgeq\lambda$.  As we are concerned with said index, however, we must build in the assumption:

\begin{thm}\label{btrsharp} Let $\lambda$ be a regular cardinal and $\cl$ an accessible category such that there exists a $\mu_\cl$-strongly compact cardinal $\kappa$.  Suppose, moreover, that $\kappa\sgeq\lambda$.  Then the powerful image of any $\lambda$-accessible functor to $\cl$ that preserves $\mu_\cl$-presentable objects is $\kappa$-accessible.\end{thm}

We note that if the $\lambda$-accessible category $\cl$ is well-$\lambda$-accessible (that is, if it is accessible in every $\kappa>\lambda$), we may dispense with the sharp inequality above.  This is significant because, as already noted, both AECs and mAECs fall under this heading.

\begin{cor}\label{btrwell} Let $\lambda$ be a regular cardinal and $\cl$ a well-accessible category such that there exists a $\mu_\cl$-strongly compact cardinal $\kappa$.  Then the powerful image of any $\lambda$-accessible functor to $\cl$ that preserves $\mu_\cl$-presentable objects is $\kappa$-accessible.\end{cor}

\section{Structural Properties Via Accessible Images}\label{sprops}

As promised in the introduction, we begin by considering the amalgamation and joint embedding properties, setting aside the disjoint versions for the moment.  Through a suitable rewriting of each one in terms of functor categories, we will move the problems to a context where we can apply Theorem~\ref{btrsharp}.  We thereby obtain Hanf numbers (in fact, the same Hanf number) for each property.

\begin{defi} Let $\ck$ be an accessible category.
\begin{enumerate}
	\item We say that $\ck$ has the {\em $<\kappa$-joint embedding property} (or {\em $<\kappa$-JEP}) if any pair of $\kappa$-presentable objects $M_0,M_1$ extends to a cospan $M_0\to N\leftarrow M_1$.  We say that $\ck$ has the {\em joint embedding property} (or {\em JEP}) if the same is true for pairs of objects of arbitrary presentability.
	\item We say that $\ck$ has the {\em $<\kappa$-amalgamation property} (or {\em $<\kappa$-AP}) if any span $N_0\stackrel{f_0}{\leftarrow} M \stackrel{f_1}{\to} N_1$ in which the objects are $\kappa$-presentable extends to a commutative square
\[
\xymatrix@=2pc{
N_0 \ar[r]^{g_0} & N \\
M \ar[u]^{f_0}\ar[r]_{f_1} & N_1\ar[u]_{g_1}
}
\]
We say that $\ck$ has the {\em amalgamation property} (or {\em AP}) if the same is true for spans of objects of arbitrary size.
\end{enumerate} 
\end{defi}

Naturally, these notions reduce to the familiar ones if $\ck$ is an AEC or mAEC.

\begin{defi}\label{catjepap} We now define the categories and functors we will use in the analysis of amalgamation and joint embedding.
	\begin{enumerate}
		\item (JEP): The category of pairs $M_0,M_1$ in $\ck$ can be identified with the category of functors from $\bullet\,\,\bullet$, the discrete category on two objects; that is, with the functor category $\ckpa$.  We may identify the category of cospans $M_0\to N\leftarrow M_1$ in $\ck$ with the category of functors from the three-object category $\bullet\to\bullet\leftarrow\bullet$ to $\ck$, i.e. $\ck^{\bullet\to\bullet\leftarrow\bullet}$.  For convenience, we will refer to the latter category as $\ck^{cosp}$.  Consider the forgetful functor
		\[F_J:\ck^{cosp}\to\ckpa\]
		that takes a cospan---a joint embedding diagram---and forgets everything but the leftmost and rightmost objects.
		
		Clearly, the image of $F_J$ is precisely the collection of jointly embeddable pairs in $\ck$.  We also note that the image of $F_J$ is closed under subobjects, and is therefore powerful.
		 		
		\item (AP): Similarly the category of spans can be identified with the functor category $\ck^{\bullet\leftarrow\bullet\to\bullet}$, which we denote by $\ck^{sp}$.  Moreover, the category of commutative squares is morally the same as the category of functors into $\ck$ from the category with four objects and morphisms arranged in a commutative square---we denote this functor category by $\ck^{csq}$.  Again, the forgetful functor
		\[F_A:\ck^{csq}\to\ck^{sp}\]
		picks out precisely the amalgamable spans.  Once again, the image of $F_A$ is itself powerful.
	\end{enumerate} 
\end{defi}

\begin{rem}\label{catchars}
(1) Notice that $\ck$ satisfies the $<\kappa$-JEP if and only if the powerful image of $F_J$ contains all pairs of $\kappa$-presentable objects, i.e. $\Preskk^{\bullet\,\bullet}$. By contrast, $\ck$ satisfies the JEP if and only if $F_J$ is an essential surjection; that is, the image is all of $\ckpa$.

(2) Similarly, $\ck$ satisfies the $<\kappa$-AP if and only if the image of $F_A$ contains $\Preskk^{sp}$. Furthermore, $\ck$ satisfies the AP if and only if $F_A$ is an essential surjection.
\end{rem}

By \cite{AR} 2.39, whenever $\ck$ is accessible and $\ca$ is a small category, $\ck^\ca$ is accessible as well, although the index of accessibility may increase.  As the categories $\ca$ invoked here are very simple, things are decidedly less complicated:

\begin{lem}\label{acclem}
	If $\ck$ is $\lambda$-accessible, so are $\ckpa$, $\ck^{sp}$, $\ck^{cosp}$, and $\ck^{csq}$.  In each case, a diagram is $\lambda$-presentable if and only if each of the objects in the diagram is $\lambda$-presentable.  Moreover, if $\ck$ is well-$\lambda$-accessible, $\ckpa$, $\ck^{sp}$, $\ck^{cosp}$, and $\ck^{csq}$ are well-$\lambda$-accessible.
\end{lem}

\proof We note that in any functor category $\ck^{\ca}$, colimits are computed componentwise, meaning that if $\ck$ has all colimits of a particular shape, so does $\ck^\ca$.  In particular, if $\ck$ is $\lambda$-accessible, it has all $\lambda$-directed colimits, and hence so does $\ck^\ca$.  So in proving $\lambda$-accessibility of our functor categories, we need only show that each contains a set of $\lambda$-presentable objects, and that any object can be realized as a $\lambda$-directed colimit of such objects.

In the case of $\ckpa$, this is clear: for any regular cardinal $\mu$, an object in $\ckpa$ is $\lambda$-presentable if and only if the pair of elements $M_0$ and $M_1$ in $\ck$ picked out by the two objects of $\bullet\,\,\bullet$ are both $\lambda$-presentable in $\ck$.  Thus there is clearly a set of $\lambda$-presentable objects, up to isomorphism.  By $\lambda$-accessiblity of $\ck$, we can express each of $M_0$ and $M_1$ as colimits of $\lambda$-directed diagrams of $\lambda$-presentable $\ck$-objects, say $M_0=\colim_{i\in I}M_0^i$ and $M_1=\colim_{j\in J}M_1^j$.  Then the diagram of pairs $(M_0^i,M_1^j)$ indexed by $I\times J$ is $\lambda$-directed, and clearly has colimit $(M_0,M_1)$.  So we are done.

For the remaining categories, we reduce the problem to the fact that if $\ck$ is $\lambda$-accessible, its category of morphisms $\ck^\to$ is also $\lambda$-accessible, and the $\lambda$-presentable morphisms are precisely those with $\lambda$-presentable domain and codomain (see \cite{AR} Ex. 2c). Indeed, since $\ck^{csq}=(\ck^\to)^\to$, we immediately get the result for the category of commutative squares. 

Let $\ck_0$ be the category obtained by formally adjoining an initial object $0$ to $\ck$. In particular, $\ck_0$ is obtained 
from $\ck$ by adding the new object $0$ and a unique morphism $0\to M$ for each $M$ in $\ck$. There are no morphisms $M\to 0$ for $M$ in $\ck$. Then $\ck_0$ is $\lambda$-accessible whenever $\ck$ is $\lambda$-accessible, and has precisely the same $\lambda$-presentable objects, with the addition of the finitely presentable object $0$. The result for $\ck^{cosp}$ follows from identifying the cospans $M_0\to N\leftarrow M_1$ in $\ck$ with the subcategory of $\ck_0^{csq}$ consisting of squares
\[
\xymatrix@=2pc{
M_0 \ar[r]^{} & N \\
0 \ar[u]^{}\ar[r]_{} & M_1\ar[u]_{}
}
\]
Working in $\ck_0^{csq}$, it is easy to verify that any cospan is a $\lambda$-directed colimit of cospans of $\lambda$-presentable objects, and that such cospans of $\lambda$-presentable objects are $\lambda$-presentable.  It follows, of course, that $\ck^{cosp}$ is $\lambda$-accessible.

Analogously, we can represent spans $N_0\leftarrow M\to N_1$ as squares
\[
\xymatrix@=2pc{
N_0 \ar[r]^{} & 1 \\
M \ar[u]^{}\ar[r]_{} & N_1\ar[u]_{}
}
\]	
in the category $\ck_1$ obtained by formally adjoining a terminal object $1$ to $\ck$. This yields the result for $\ck^{sp}$.
	
	The moreover clause follows in a similar fashion.
\qed

From this perspective, the proofs of the existence of a Hanf number for JEP and AP in accessible categories are essentially trivial, modulo Theorem~\ref{btrsharp}.

\begin{thm}\label{thmJEP} Let $\ck$ be a $\lambda$-accessible category, and let $\kappa$ be a $\mu_\ck$-strongly compact cardinal with $\kappa\sgeq\lambda$.  If $\ck$ has the $<\kappa$-JEP, then it has the JEP.\end{thm}
\proof
	Consider the functor $F_J:\ck^{cosp}\to\ckpa$.  Each of $\ckpa$ and $\ck^{cosp}$ are $\lambda$-accessible, by Lemma~\ref{acclem}, and clearly $F_J$ preserves $\lambda$-directed colimits.  That is, $F_J$ is $\lambda$-accessible.  
	
	Notice that, since $\ck$ and $\ckpa$ are both $\lambda$-accessible and have subcategories of $\lambda$-presentables of the same size, $\mu_\ck=\mu_{\ckpa}$.  As the $\mu_\ck$-presentable objects in $\ck^{cosp}$ are precisely the cospans whose objects are $\mu_\ck$-presentable, they are certainly preserved by $F_J$.
	
	Hence by Theorem~\ref{btrsharp}, the powerful image of $F_J$ is $\kappa$-accessible.  As noted in Definition \ref{catjepap}(1), this powerful image consists of precisely the jointly embeddable pairs of objects in $\ck$.  Let $(M_0,M_1)$ be an arbitrary pair of objects in $\ck$, i.e. an arbitrary object of $\ckpa$.  Because $\kappa\sgeq\lambda$, $\ckpa$ is $\kappa$-accessible, meaning that $(M_0,M_1)$ in $\ckpa$ can be realized as a $\kappa$-directed colimit of pairs of $\kappa$-presentable objects; that is,
	\[(M_0,M_1)=\colim_{i\in I}(M_0^{i},M_1^{i})\]
	with $I$ a $\kappa$-directed set.  
	
	Since $\ck$ satisfies the $(<\kappa)$-JEP, by assumption, Remark \ref{catchars}(1) implies that each $(M_0^{i},M_1^{i})$ lies in the powerful image of $F_J$.  As the powerful image of $F_J$ is $\kappa$-accessible, it is closed under $\kappa$-directed colimits, and therefore contains $(M_0,M_1)$.  In particular, $M_0$ and $M_1$ admit a joint embedding.  So $\ck$ satisfies the JEP.\qed

\begin{thm}\label{thmAP} Let $\ck$ be a $\lambda$-accessible category, and let $\kappa$ be a $\mu_\ck$-strongly compact cardinal with $\kappa\sgeq\lambda$.  If $\ck$ has the $<\kappa$-AP, then it has the AP.\end{thm}
\proof The theorem follows by the same argument as Theorem~\ref{thmJEP}, with $\ck^{sp}$ and $\ck^{csq}$ in place of $\ckpa$ and $\ck^{cosp}$, respectively.\qed

We have simpler statements in case $\ck$ is well-accessible:

\begin{thm}\label{thmJEPAPwell}
	Let $\ck$ be a well $\lambda$-accessible category, and let $\kappa$ be a $\mu_\ck$-strongly compact cardinal.  If $\ck$ has the $<\kappa$-JEP ($<\kappa$-AP), then it satisfies the JEP (AP).
\end{thm}
\proof This follows by Corollary~\ref{btrwell} and the moreover clause of Lemma~\ref{acclem}.\qed

\begin{rem}\label{BalBoncomp}
	We note that Theorem 2.0.5 in \cite{BaBo}, which covers the case in which $\ck$ is an AEC, assumes that $\kappa$ is a strongly compact cardinal greater than $\lsk$.  In particular, $\kappa$ must be strongly inaccessible.  In this case $\ck$ is $\lsk^+$-accessible, and the cardinals $\lambda=\lsk^+$, $\beta=|\Preslk|$, $\gamma_\ck=\max(\lambda,\beta)$, and $\mu_\ck=(\gamma_\ck^{<\gamma_\ck})^+$ are all less than $\kappa$.  It follows that $\kappa$ is $\mu_\ck$-strongly compact.  Thus the specialization of Theorem~\ref{thmJEPAPwell} to the AEC case involves, {\em a priori}, a weaker large cardinal assumption: we note, however, that, per \cite{BU}~2.2, it remains open whether this assumption is genuinely (consistently) weaker.
\end{rem}

\section{Disjoint Analogues}\label{disj}

Capturing the disjoint versions of amalgamation and joint embedding requires more care, insofar as we must frame the notion of disjointness in an abstract category.  The disjoint JEP for an AEC $\ck$, for example, asserts that any models $M_0$ and $M_1$ can be completed to a cospan $M_0\stackrel{g_0}{\to} N\stackrel{g_1}{\leftarrow} M_1$ where the ranges of the embeddings are disjoint, i.e. the following diagram is a pullback in $\Set$:
\[
\xymatrix@=2pc{
UM_0 \ar[r]^{U(g_0)} & UN \\
\emptyset \ar[u] \ar[r]  & UM_1\ar[u]_{U(g_1)}
}
\]
Here the essential role is played by the empty set, which is a {\em strictly initial object} in $\Set$: it is initial, and any morphism $M\to\emptyset$ is an isomorphism. The difficulty is that we cannot expect $\ck$ itself to have a strict initial object (or, indeed, pullbacks), so we cannot necessarily pull this diagram back to one in $\ck$, considered as an abstract category.  

We skirt this difficulty by broadening the notion of disjointness to incorporate categories that may lack a strict initial object: A cospan $M_0\to N\leftarrow M_1$ in a category $\ck$ will be called \textit{disjoint} if it cannot be completed to a commutative square or
\[
\xymatrix@=2pc{
M_0 \ar[r] & N \\
\emptyset \ar[u] \ar[r]  & M_1\ar[u]
}
\]
is the only completion, where $\emptyset$ is strictly initial in $\ck$. Note that this square, if it exists, is a pullback square in $\ck$. 

The category of disjoint cospans in $\ck$ will be denoted 
$\ck^{dcosp}$, while $\ck^{psq}$ will be the category of pullback squares in $\ck$.

\begin{lem}\label{disjoint1}
Let $\ck$ be a $\lambda$-accessible category. Then the category $\ck^{dcosp}$ is $\lambda$-accessible and $\ck^{psq}$ is
$\lambda^+$-accessible.  
\end{lem}
\proof
Any morphism in $\ck^{cosp}$ to a disjoint cospan has disjoint domain. Thus the result for $\ck^{dcosp}$ follows from
that for $\ck^{cosp}$ (see \ref{acclem}). 

Following \cite{MP} 6.2.1(i), $\ck^{psq}$ is given by the pseudopullback
\[
\xymatrix@=2pc{
\ck^{cosp} \ar[r]^{H} & \Set^{\ca^{\op}} \\
\ck^{psq} \ar[u]^{}\ar[r]_{} & \ck\ar[u]_{E}
}
\]
where $E$ is the canonical functor $\ck\to\Set^{\ca^{\op}}$ given by $EM=\ck(-,M)\upharpoonright\ca^{\op}$, where $\ca$ denotes $\Preslk$ (see \cite{AR} 1.25), and $H$ is the composition
$$
\ck^{cosp}\to (\Set^{\ca^{\op}})^{cosp}\to\Set^{\ca^{\op}}
$$
where the first functor is $E^{sp}$ and the second functor sends a cospan in the presheaf category $\Set^{\ca^{\op}}$ to its pullback---note that, while $\ck$ may not have pullbacks, $\Set^{\ca^{\op}}$ always will. Since $E$ and $H$ preserve $\lambda$-directed
colimits and $\lambda$-presentable objects, $\ck^{psq}$ is $\lambda^+$-accessible (see \cite{RR} 2.2).
%
\qed


\begin{defi}\label{disjoint2}
A category $\ck$ has the \textit{disjoint embedding property} (or {\em DJEP}) if any pair of objects $M_0$ and $M_1$ can be completed
to a disjoint cospan $M_0\to N\leftarrow M_1$.

$\ck$ has the \textit{disjoint amalgamation property} (or {\em DAP}) if any span $N_0\leftarrow M\rightarrow N_1$ can be completed to a pullback
square
\[
\xymatrix@=2pc{
N_0 \ar[r] & N \\
M \ar[u] \ar[r]  & N_1\ar[u]
}
\]
The meanings of {\em $<\kappa$-$DJEP$} and {\em $<\kappa$-$DAP$} are clear.\end{defi} 

\begin{rem}\label{disjgens} Those familiar with the AEC versions of these properties will recognize our formulation of the DAP as a straightforward generalization.  The DJEP is as formulated here is slightly more delicate: in fact, even in case $\ck$ is an AEC, it yields a notion still more general than the {\it nearly disjoint JEP}, in which any two models can be embedded in such a way that they overlap only on the submodel generated from the constant symbols.  It suffices for our purposes here, however.  That is, if an AEC $\ck$ satisfies the $<\kappa$-DJEP for any $\kappa$, then it cannot contain any strictly initial object other than the empty model.  If it contains no strictly initial object, we formally adjoin the empty model to serve this purpose, a step which does not change the accessibility of $\ck$.  In either case, we have the $<\kappa$-DJEP in the sense of Definition~\ref{disjoint2}.  When we show below that $\ck$ therefore satisfies the global DJEP, again in the sense of Definition~\ref{disjoint2}, we can be certain that this translates into the familiar DJEP for $\ck$ as an AEC, simply because we still have the empty model as the strictly initial object.\end{rem}

\begin{thm}\label{disjoint3}
Let $\ck$ be a $\lambda$-accessible category. If $\kappa$ is a $\mu_{\ck,\lambda}$-strongly compact cardinal with $\kappa\sgeq\lambda$ and $\ck$ has the $<\kappa$-DJEP, then it satisfies the DJEP. If $\kappa$ is $\mu_{\ck,\lambda^+}$-strongly compact with $\kappa\sgeq\lambda^+$ and $\ck$ has 
the $<\kappa$-DAP, then it satisfies the DAP.
\end{thm}
\proof
For DJEP we use the $\lambda$-accessible functor
\[\ck^{dcosp}\to\ckpa\]
in the same way as in \ref{thmJEP}. For DAP we use the functor \[\ck^{psq}\to\ck^{sp}\] from the proof
of \ref{disjoint1}. Since this functor is $\lambda^+$-accessible (see \cite{CR} 3.1), we proceed as in \ref{thmAP}.\qed

In the case of DAP, $\mu_\ck$ is taken for $\lambda^+$ and not for $\lambda$. 
Again, we have simpler statements for well-accessible categories.

\begin{thm}\label{disjoint4} 
Let $\ck$ be a well $\lambda$-accessible category, and let $\kappa$ be $\mu_{\ck,\lambda}$-strongly compact.  If $\ck$ has the $<\kappa$-DJEP, then it satisfies the DJEP.  Similarly, if $\kappa$ is $\mu_{\ck,\lambda^+}$-strongly compact and $\ck$ has the $<\kappa$-DAP, it has the AP.
\end{thm}

\begin{rem}
	As noted in Remark~\ref{BalBoncomp}, the specialization of Theorem~\ref{disjoint4} to AECs involves an {\em a priori} weaker large cardinal hypothesis than \cite{BaBo} 2.0.5.
\end{rem}

\section{Metric Tameness}\label{mtame}

We turn now to the most technical and involved of the results to be presented here, the $d$-tameness of metric AECs under the assumption of sufficiently strongly compact cardinals.  We will again provide the necessary basic definitions; extensive background and motivation for this question, and for the broader study of mAECs, can be found in, e.g. \cite{HH} and \cite{VZ}.

As suggested in the introduction, the method to be used here is a generalization of that employed in \cite{LR} to establish the discrete version of this result, namely the tameness of AECs under the assumption of arbitrarily large strongly compact cardinals.  We dwell briefly on that case here, largely as an opportunity to recast it in a form better suited to the needed generalization.  In the discrete context, a Galois type over a model $M$ is an equivalence class of pairs $(f,a)$, where $f:M\to N$ is a $\ck$-embedding and $a\in UN$.  Here two pairs $(f_0,a_0)$ and $(f_1,a_1)$ are said to be equivalent if there is an amalgam of $f_0$ and $f_1$ in which $a_0$ and $a_1$ are identified, i.e. there are $g_0:N_0\to N$ and $g_1:N_1\to N$ such that the following diagram commutes:

\[\xymatrix@=.5pc@C=.1pc{
a_0\ar@{|->}[rrrr] & & & & U(g_0)(a_0)\ar@{=}[dr] & \\
UN_0\ar[rrrr]_{U(g_0)} & & & & UN & {U(g_1)(a_1)}\\
N_0\ar@{}[rrrddd]|{\bigcirc}\ar@{-->}[rrr]_{g_0} & \hspace{.9 cm} & & N & & \\
& & & & & \\
& & & & & \\
M\ar[uuu]_{f_0}\ar[rrr]^{f_1} & & & N_1\ar@{-->}[uuu]^{g_1} & UN_1\ar[uuuu]_{U(g_1)} & a_1\ar@{|->}[uuuu]
}\]

This is an equivalence relation provided $\ck$ has amalgamation.  Equivalently, two pairs $(f_0,a_0)$ and $(f_1,a_1)$ represent distinct types just in case for any amalgam $g_0:N_0\to N$ and $g_1:N_1\to N$, there is a monomorphism $h$ from the two point set $2=\{0,1\}$ to $UN$ such $h(0)=U(g_0)(a_0)$ and $h(1)=U(g_1)(a_1)$, as in the diagram below.

\[\xymatrix@=.6pc@C=.3pc{
& & & & & & 0\ar@{|->}[ddll]\hspace{.5 cm} & \vspace{.3 cm}\\
& & & & & & \hspace{.1 cm}2\ar@{-->}[ddll]_h\hspace{-.55 cm} & \vspace{-.2 cm}\\
a_0\ar@{|->}[rrrr] & & & & U(g_0)(a_0) & & & {1}\ar@{|->}[dll]+UR\\
UN_0\ar[rrrr]_{U(g_0)} & & & & UN & {U(g_1)(a_1)} & & \\
N_0\ar@{}[rrrddd]|{\bigcirc}\ar[rrr]_{g_0} & \hspace{.9 cm} & & N & & & & \\
& & & & & & & \\
& & & & & & & \\
M\ar[uuu]_{f_0}\ar[rrr]^{f_1} & & & N_1\ar[uuu]^{g_1} & UN_1\ar[uuuu]_{U(g_1)} & a_1\ar@{|->}[uuuu] & & 
}\]

Although the argument of \cite{LR} is built around the former characterization of equivalence, the latter is of greater significance for the discussion here: we measure distinctness of types in AECs by monomorphisms into amalgams from $2$, the natural source for such probes in $\Set$.  What if, instead, we consider the metric case, where structures have underlying complete metric spaces, rather than sets?  We will want to measure not inequivalence of pairs, but the distance between them.  Hence the natural objects with which to probe in $\Met$ are the two point metric spaces of diameter $\epsilon$, $2_{\epsilon}$, for $\epsilon\in\reals_{>0}$.

There is a great deal to be made precise.  We recall first that for all the minor differences between metric and discrete AECs---chiefly that the operative notion of size is density character rather than cardinality and that directed colimits need not be concrete---the fundamental definition of Galois types is precisely the same: as equivalence classes of pairs $(f_0:M\to N_0,a_0)$, $(f_1:M\to N_1,a_1)$.  One can, however, speak not merely of the equality or inequality of types, as in the discrete case, but of the distance between them.  This distance is typically defined as the Hausdorff distance between the corresponding orbits in a monster model $\monst$: the assumption of AP, JEP, and the existence of arbitrarily large models guarantees that $\monst$ exists; the {\em perturbation property} (\cite{HH}), also known as the {\em continuity of types property} (\cite{VZ}), guarantees that this distance is a metric.  We make precisely the same assumptions here:

\begin{asm}\label{maecassump} From now on, we assume that our mAEC $\ck$ has the AP and JEP, has arbitrarily large models, and satisfies the perturbation property.\end{asm}

We need an alternative, albeit equivalent, definition of the distance between types, however, for reasons which will soon become clear.

\begin{defi} Let $\ck$ be an mAEC (satisfying the conditions of Assumption~\ref{maecassump}).  Given types $p_0=\gatp(f_0,a_0)$ and $p_1=\gatp(f_1,a_1)$ over $M\in\ck$, $f_i:M\to N_i$, we define the {\em distance between $p_0$ and $p_1$}, denoted $d(p_0,p_1)$, to be the infimum of the distances between the images of $a_0$ and $a_1$ in all amalgams of $f_0$ and $f_1$.  That is,
\[d(p_0,p_1)=\inf\{d_N(U(g_0)(a_0),U(g_1)(a_1))\,|\,g_i:N_i\to N\hbox{ and }g_0f_0=g_1f_1\}\]
\end{defi}

Notice that $d_N(U(g_0)(a_0),U(g_1)(a_1))=\epsilon$ if and only if there is an isometry $h:2_\epsilon\to N$, where 
$2_\epsilon=\{0,1\}$ is the two point space of diameter $\epsilon$, such that $h(0)=U(g_0)(a_0)$ and $h(1)=U(g_1)(a_1)$.  
This is witnessed by a diagram of the following shape:

\[\xymatrix@=.6pc@C=.3pc{
& & & & & & 0\ar@{|->}[ddll]\hspace{.5 cm} & \vspace{.3 cm}\\
& & & & & & \hspace{.1 cm}2_{\epsilon}\ar@{-->}[ddll]_h\hspace{-.55 cm} & \vspace{-.2 cm}\\
a_0\ar@{|->}[rrrr] & & & & U(g_0)(a_0) & & & {1}\ar@{|->}[dll]+UR\\
UN_0\ar[rrrr]_{U(g_0)} & & & & UN & {U(g_1)(a_1)} & & \\
N_0\ar[rrr]_{g_0} & \hspace{.9 cm} & & N & & & & \\
& & & & & & & \\
& & & & & & & \\
M\ar[uuu]_{f_0}\ar[rrr]^{f_1} & & & N_1\ar[uuu]^{g_1} & UN_1\ar[uuuu]_{U(g_1)} & a_1\ar@{|->}[uuuu] & & 
}\]

\begin{defi}\label{dtamecats} For each $\epsilon>0$, we define $\cl_\epsilon$ to be the category of (pointed) commutative squares in $\Met$ equipped with a specified such isometry $h$.  Let $\cl$ denote the category of pairs $(f_0,a_0)$ and $(f_1,a_1)$ or, in the notation of \cite{LR}, $(f_0,f_1,a_0,a_1)$.  Finally, we denote by $G_\epsilon:\cl_\epsilon\to\cl$ the obvious forgetful functor.\end{defi}  

In the proof of $d$-tameness that follows, we will use the eventual $\kappa$-accessibility of the powerful images of these functors $G_\epsilon$ to show that if two types over a model are separated by distance at least $\epsilon$, then this fact is witnessed by restrictions to a model of density character $\kappa$, in the sense that the restrictions are also at distance at least $\epsilon$.  That is, we show that any mAEC is strongly $d$-tame:

\begin{defi} We say that an mAEC $\ck$ is {\em $\kappa$-$d$-tame} if for every $\epsilon>0$ there is a $\delta>0$ such that for any $M\in\ck$ and Galois types $p,q$ over $M$, if $d(p,q)\geq\epsilon$ then there exists $N\ksst M$ with dc$(N)\leq\kappa$ such that $d(p\upharpoonright N,q\upharpoonright N)\geq\delta$.  We say that $\ck$ is {\em strongly $\kappa$-$d$-tame} if the above holds with $\delta=\epsilon$.  

We say that $\ck$ is {\em (strongly) $d$-tame} if it is (strongly) $\kappa$-$d$-tame for some $\kappa$.\end{defi}

\begin{thm} Suppose that for every cardinal $\mu$ there exists a $\mu$-strongly compact cardinal. Then every mAEC is strongly $d$-tame.\end{thm}

\proof Let $\ck$ be an mAEC, and let $U:\ck\to\Met$ be the forgetful functor to the category of complete metric
spaces and isometric embeddings. Let $\cl$, $\cl_\epsilon$, and $G_\epsilon:\cl_\epsilon\to\cl$ be as described in Definition~\ref{dtamecats}.  Then both $\cl,\cl_\varepsilon$ are accessible and the forgetful functor $G_\varepsilon:\cl_\varepsilon\to\cl$ is accessible as well. In fact, $\cl$ is the category of commutative squares 
\[
\xymatrix@=2pc{
N_0 \ar[r]^{g_0} & N \\
M \ar[u]^{f_0}\ar[r]_{f_1} & N_1\ar[u]_{g_1}
}
\]
equipped with elements $a_0\in UN_0$ and $a_1\in UN_1$. Thus its accessibility follows from that of $\ck^{csq}$ and from
the preservation of $\aleph_1$-directed colimits by the forgetful functor $U:\ck\to\Met$. Similarly for $\cl_\epsilon$ where
one adds the condition that $d_N(g_0(a_0),g_1(a_1))=\epsilon$.
Clearly, if the quadruple $(f_0,f_1,a_0,a_1)$ belongs to the image of $G_\varepsilon$ then 
the distance of the corresponding Galois types is $\leq\varepsilon$. Thus the distance $d(p_0,p_1)$ of Galois types $p_0=(f_0,a_0)$ and $p_1=(f_1,a_1)$ is equal to
the infimum of the set 
\[\{\epsilon>0\,|\,(f_0,f_1,a_0,a_1)\in G_\varepsilon(\cl_\varepsilon)\}\]

It is easy to see that the full image of $G_\varepsilon$ is a sieve: for a morphism $(u,v,w):(f^{'}_0,f^{'}_1,a_0',a_1')\to (f_0,f_1,a_0,a_1)$, i.e. 
\[
\xymatrix@=2pc{
UN_0'\ni a_0'\ar@{|->}@<3.3ex>[d]_{U(u)} & N_0' \ar[d]_u & M\ar[r]^{f_1'}\ar[l]_{f_0'}\ar[d]_v & N_1'\ar[d]_w & a_1'\in UN_1'\ar@{|->}@<-4ex>[d]_{U(w)}\\
UN_0\ni a_0 & N_0  & M\ar[r]_{f_1}\ar[l]^{f_0} & N_1 & a_1\in UN_1\\
}
\]
with $(f_0,f_1,a_0,a_1)\in G_\varepsilon(\cl_\varepsilon)$, we have $(f^{'}_0,f^{'}_1,a_0',a_1')\in G_\varepsilon(\cl_\varepsilon)$ as well. Thus the full image of $G_\varepsilon$ is powerful. Let $\lambda$ 
be a regular cardinal such that all functors $G_\varepsilon$ are $\lambda$-accessible and preserve $\lambda$-presentable objects. Let $\mu_\cl$ be the cardinal from defined in Section~\ref{prelim}---just before Theorem~\ref{btrsharp}---and let $\kappa$ be a $\mu_\cl$-strongly compact cardinal. Then, following Theorem~\ref{btrsharp}, the full image of any $G_\varepsilon$ is 
$\kappa$-accessible. 

Assume that the distance $d(p_0\chi,p_1\chi)<\varepsilon$ for each morphism $\chi:X\to M$ in $\ck$ with $X$ $\kappa$-presentable, where we write $p_i\chi$ as an abbreviation for the type represented by $(f_i\chi,a_i)$. We claim that there are cofinally many $\chi:X\to M$ giving the same distance $d(p_0\chi,p_1\chi)=\delta$ for some $\delta<\varepsilon$. In fact, assuming the contrary, for any $\delta<\varepsilon$ there is $\chi_\delta:X_\delta\to M$
with no factorization $\chi_\delta=\chi't_\delta$ such that $d(p_0\chi',p_1\chi')=\delta$. There are at most $2^{\aleph_0}$ many witnesses $\chi_\delta$, and given that $M$ is a $\kappa$-directed colimit of the $\chi:X\to M$, there must therefore be a map $\chi':X'\to M$ through which all of the $\chi_\delta$ factor.  That is, for every $\delta<\varepsilon$, there is $t_\delta:X_\delta\to X'$ such that the following diagram commutes:
\[
\xymatrix@=4pc{
X'\ar[r]^{\chi'} & M\\
X_\delta\ar[u]^{t_\delta}\ar[ur]_{\chi_\delta}
}
\]
Note that $d(p_0\chi',p_1\chi')=\delta'$ for some $\delta'<\varepsilon$.  Of course, $\chi_{\delta'}$ must factor through $\chi'$, which is a contradiction.

Let $\cd$ denote the resulting cofinal set. For each $\chi\in\cd$ , there is a sequence of $\delta_n>\delta$ converging to $\delta$ such that 
$(p_0\chi,p_1\chi)$ belongs to the image of $G_{\delta_n}$ for each $n$. Since $(2^{\aleph_0})^{\aleph_0}=2^{\aleph_0}<\kappa$, there is a cofinal subset of $\cd$ realizing the same sequence $\delta_n$. Otherwise, for any sequence $\delta_\ast=(\delta_n)_n$ there is $\chi_{\delta_\ast}$ with no factorization 
$\chi_{\delta_\ast}=\chi't$ such that $(p_0\chi',p_1\chi')$ does not belong to the image of $G_{\delta_n}$ for each $n$. Since all of the $\chi_{\delta_\ast}$ factorize through some $\chi'$ we get a contradiction, as above. Thus $(p_0,p_1)$ belongs to the image of $G_{\delta_n}$ for each $n$, by accessibility of the images. Hence 
$d(p_0,p_1)\leq\delta$. Consequently, $d(p_0,p_1)<\varepsilon$.\qed

Noting that a strongly $d$-tame $\ck$ is $d$-tame in the sense of \cite{Z}, we have:

\begin{cor}
Suppose that for every cardinal $\mu$ there exists a $\mu$-strongly compact cardinal. Then every mAEC is $d$-tame.
\end{cor}

This last result forms half of Theorem~3.11 in \cite{BoZ}, although it appears that their proof should entail strong $d$-tameness as well.  A more substantive difference is that, while they rely on a delicate argument involving metric ultraproducts, the argument above more clearly illustrates the common underlying structure of the discrete and continuous problems, and points the way to further generalizations.

\section{Extensions}

We note two clear directions in which the results of this paper can be extended.  First, the focus on joint embedding and amalgamation in Sections~\ref{sprops} and \ref{disj} is due entirely to practical considerations: these properties, amounting to the existence of coproduct-like cones over pairs and pushout-like cones over spans, are the ones of greatest interest to model theorists.  The form of argument via accessible images is vastly more general, though.  For example, while model theorists work almost exclusively with monomorphisms, meaning that coequalizer-like (or quotient-like) cones over pairs of parallel arrows are rarely of interest, these make perfect sense in a general accessible category.  We might say, for example, that an accessible category has the coequalizing property (CP) if every parallel pair of arrows $f,g:M\to N$ is coequalized by a map $h:N\to N'$, i.e. $hf=hg$.  The existence of such an $h$ can again be formulated in terms of the image of a forgetful functor $F_C:\ck^{\ca'}\to\ck^\ca$, with $\ca'=\bullet\rightrightarrows\bullet\to\bullet$ (with both compositions equal) and $\ca=\bullet\rightrightarrows\bullet$.  As the image is again closed under subobjects, i.e. powerful, we obtain a Hanf number for the CP, just as in Theorem~\ref{thmJEP} and ~\ref{thmAP}.  There are is an abundance of possible analogues, in fact, limited only by this last stipulation: the image of the forgetful functor, i.e. the class of completable diagrams, must be closed under subobjects.

We note also that ongoing work with Zambrano suggests that the results of Section~\ref{mtame} can be extended well beyond the case of mAECs.  Of particular interest are classes of structures with underlying $Q$-sets, with $Q$ a suitably cocontinuous quantale, in roughly the sense of \cite{F}.  This would mean immediate tameness results across a wide array of fields, including fuzzy semantics, model theory over probability spaces, and model theory in categories of sheaves over sufficiently well-behaved spaces.


\begin{thebibliography}{10}
 
\bibitem{AR} J. Ad\'{a}mek and J. Rosick\'{y}. {\em Locally Presentable and Accessible Categories}. Cambridge University
Press, 1994.
 
\bibitem{Ba} J. Baldwin. {\em Categoricity}. AMS, 2009.

\bibitem{BaBo} J. Baldwin and W. Boney. Hanf numbers and presentation theorems in AECs.  To appear in {\em Beyond First Order Model Theory}, ed. J. Iovino, CRC Press.  arXiv:1511.02935v2.

\bibitem{BR} T. Beke and J. Rosick\'y. Abstract elementary classes and accessible categories. {\em Ann. Pure Appl. Log.}, 163:2008-2017, 2012.

\bibitem{B} W. Boney. Tameness from large cardinal axioms. {\em J. Symb. Log.}, 79(4):1092-1119, 2014.

\bibitem{muaecs} W. Boney, R. Grossberg, M. Lieberman, J. Rosick\'y, and S. Vasey. $\mu$-abstract elementary classes and other generalizations. {\em J. Pure Appl. Alg.}, 220(9):3048-3066, 2016.

\bibitem{BU} W. Boney and S. Unger. Large cardinal axioms from tameness in AECs. To appear in {\em Proc. AMS}.  arXiv:1509.01191v3.

\bibitem{BoZ} W. Boney and P. Zambrano. Around the set-theoretical consistency of $d$-tameness of metric abstract elementary classes. arXiv:1509.01191v1.

\bibitem{BTR} A. Brooke-Taylor and J. Rosick\'y. Accessible images revisited. {\em Proc. AMS}, 145(3):1317-1327.

\bibitem{CR} B. Chorny and J. Rosick\'y. Locally class-presentable and class-accessible categories. {\em J. Pure Appl. Alg.}, 216:
2113-2125, 2012.
                   
\bibitem{EM} P. Eklof and A. Mekkler. {\em Almost Free Modules: Set-theoretic methods}. North Holland, 1990.

\bibitem{HH} A. Hirvonen and T. Hyttinen. Categoricity in homogeneous complete metric spaces. {\em Arch. Math. Log.}, 48:269-322, 2009.

\bibitem{F} R. Flagg. Quantales and continuity spaces. {\em Alg. Univ.}, 37:257-276, 97.

\bibitem{L} M. Lieberman. Category theoretic aspects of abstract elementary classes. {\em Ann. Pure Appl. Log.}, 162:903-915, 2011.

\bibitem{LR} M. Lieberman and J. Rosick\'y. Classification theory for accessible categories. {\em J. Symb. Log.}, 81(1):151-165, 2016.

\bibitem{LRmetr} M. Lieberman and J. Rosick\'y. Metric abstract elementary classes as accessible categories.  To appear in the {\em J. Symb. Log.}.  arXiv:1504.02660v5.
 
\bibitem{MP} M. Makkai and R. Par\' e. {\em Accessible Categories: The Foundations of Categorical Model Theory}. AMS, 1989.

\bibitem{RR} G. Raptis and J. Rosick\'y. The accessibility rank of weak equivalences. {\em Theory Appl. Categ.}, 30:687-703, 2015.    

\bibitem{R} J. Rosick\'y. Accessible categories, saturation and categoricity. {\em J. Symb. Log.}, 62:891-901, 1997.

\bibitem{VZ} A. Villaveces and P. Zambrano. Limit models in metric abstract elementary classes: the categorical case. {\em Math. Log. Quart.}, 62(4-5):319-334, 2016.

\bibitem{Z} P. Zambrano. Around superstability in metric abstract elementary classes. PhD thesis.


\end{thebibliography}
\end{document}